\documentclass[preprints,article,accept,pdftex,moreauthors]{Definitions/mdpi} 
\firstpage{1} 
\pubvolume{1}
\issuenum{1}
\articlenumber{0}
\pubyear{2022}
\copyrightyear{2022}
\datereceived{} 
\dateaccepted{} 
\datepublished{} 
\hreflink{https://doi.org/} 

\RequirePackage{algorithm}
\RequirePackage{algorithmic}
\RequirePackage{mathtools}
\DeclareMathOperator{\Null}     {Null}
\DeclareMathOperator{\poly}     {poly}
\DeclareMathOperator{\rank}     {rank}
\DeclareMathOperator{\trace}    {tr}

\Title{An Inexact Feasible Quantum Interior Point Method for Linearly Constrained Quadratic Optimization}

\TitleCitation{}

\Author{Zeguan Wu $^{1}$\orcidA{}, Mohammadhossein Mohammadisiahroudi $^{1}$\orcidB{}, Brandon Augustino $^{1}$\orcidC{}, Xiu Yang $^{1}$\orcidD{} and Tam\'as  Terlaky $^{1,}$*\orcidE{}}

\longauthorlist{yes}

\AuthorNames{}

\AuthorCitation{}

\address[1]{%
$^{1}$ \quad Department of Industrial and Systems Engineering, Lehigh University}

\corres{Correspondence: terlaky@lehigh.edu}

\abstract{Quantum linear system algorithms (QLSAs) have the potential to speed up algorithms that rely on solving linear systems. Interior Point Methods (IPMs) yield a fundamental family of polynomial-time algorithms for solving optimization problems. IPMs solve a Newton linear system at each iteration to find the search direction, and thus QLSAs can potentially speed up IPMs. Due to the noise in contemporary quantum computers, such quantum-assisted IPM (QIPM) only allows an inexact solution for the Newton linear system. Typically, an inexact search direction leads to an infeasible solution. In our work, we propose an Inexact-Feasible QIPM (IF-QIPM) and show its advantage in solving linearly constrained quadratic optimization problems. We also apply the algorithm to $\ell_1$-norm soft margin support vector machine (SVM) problems and obtain the best complexity regarding dependence on dimension. This complexity bound is better than any existing classical or quantum algorithm that produces a classical solution.}

\keyword{Quantum Computing; Interior Point Method; Quadratic Optimization}

\MSC{90C20; 90C51; 81P68}
\begin{document}




\section{Introduction}
Linearly constrained quadratic optimization (LCQO) is defined as optimizing a convex quadratic objective function over a set of linear constraints. This problem reduces to linear optimization when the the quadratic objective function is linear. LCQO has rich theory, algorithms, and applications. Many machine learning problems are LCQO problems, including variants of least square problems and variants of support vector machine problems \cite{nocedal1999numerical, boser1992training}. Some important optimization algorithms also have LCQO subproblems, e.g. sequential quadratic programming \cite{nocedal1999numerical}. 

 The modern age of IPMs launched by Karmarkar's invention of the projective method for linear optimization (LO). Since then, a lot of variants of IPMs have been studied for not only LO problems but also for nonlinear optimization problems, including LCQO problems \cite{roos1997theory,polik2010interior}.

Contemporary IPMs look for the optimal solution by moving in a neighbourhood of the central path. IPMs can be divided into two classes: feasible or infeasible. Feasible IPMs start with a feasible solution and keep feasibility; infeasible IPMs start with an infeasible interior solution and so do not require a feasible solution to start with. For LCQO problems with $n$ variables, some feasible IPMs produce an $\epsilon$-approximate solution in at most $\mathcal{O}(\sqrt{n}\log(1/\epsilon))$ IPM iterations, while infeasible IPMs require $\mathcal{O}(n^2\log(1/\epsilon))$ IPM iterations to generate an $\epsilon$-approximate solution \cite{gondzio2013convergence,lu2006iterative}.

At each IPM iteration a linear system needs to be solved to obtain the search direction, called the Newton direction. Such a Newton linear system is traditionally in the form of augmented system or the normal equation system. Classically these linear systems can be solved exactly using Bunch-Parlett factoriztion if the matrices in the systems are symmetric indefinite \cite{bunch1971direct}, or Cholesky factorization if the matrices are symmetric positive definite. The complexity of solving the linear systems is $\mathcal{O}(n^3)$. The linear systems can also be solved inexactly using some inexact methods, e.g., Krylov subspace methods. Such inexact methods might take less iterations if the desired accuracy of the solutions to the linear systems is not high. But such inaccuracy of the solutions to the linear systems, i.e., inaccuracy of the search directions, might result in  infeasibility of the solutions generated by IPMs. To maintain feasibility of solutions,  \cite{mohammadisiahroudi2021inexact} introduces the so-called orthogonal subspace system (OSS) for LO problems. A feasible solution can be recovered from an inexact solution to OSS. We extend their OSS for LO prolems to LCQO problems and provide an efficient method to construct the OSS. With the OSS, we can obtain an inexact feasible IPM -- solving for search direction inexactly but maintaining the feasibility of solution throughout the process of our IPM. The feasibility of solution gives better IPM iteration complexity and the bottleneck becomes solving the linear system, OSS.

With the development of quantum technology, many quantum-assisted algorithms have been proposed for many optimization problems. Following the invention of quantum algorithms for solving linear systems of equations \cite{harrow2009quantum}, many researchers are encouraged to study whether QLSAs would yield quantum speedups in classical algorithms. In particular, QIPMs have been proposed for for LO problems \cite{kerenidis2020quantum, mohammadisiahroudi2022efficient} and semidefinite optimization problems \cite{augustino2021quantum} that utilize QLSAs to solve the Newton linear system that arises in each iteration of IPMs. Similar ideas have also been applied to accelerate the solution of some machine learning applications, such as linear regression \cite{SchuldRegression} and the support vector machine training problem \cite{kerenidis2021quantum}. However, linearly constrained quadratic optimization problems, which are fundamental to both optimization and machine learning, have not been formally studied in the quantum literature yet. 

The remaining part of this paper is organized as follow: in Section 2, we introduce IPMs for LCQO and the OSS system; in Section 3, we discuss how to use quantum algorithms to find the Newton directions and analyze the complexity of our IF-QIPM; in Section 4, we apply our IF-QIPM to support vector machine problem. Discussions are provided in Section 5, and some technical proofs are moved to the Appendix.

\section{Preliminary}
\subsection{Notations}
In this section, we introduce notations we use.
Vectors are typically represented by lower-case letters. For $n$-dimensional all-zero vector, we represent it with $0_{n}$ if the dimension is $n$, or simply $0$ if the dimension is obvious in the context. For $n$-dimensional all-one vector, we represent it with $e_n$, or simply $e$ if the dimension is obvious in the context.

Matrix are typically represented with upper-case letters. For $n$-dimensional identity matrix, we represent it with $I_{n\times n}$, or simply $I$ if the dimension is obvious in the context. For $n\times m$-dimensional all-zero matrix, we represent it with $0_{n\times m}$, or simply $0$ if the dimension if obvious in the context. For a general $n\times m$-dimensional matrix $H$, we represent its $i$th row by $H_{i\cdot}$ and $j$th column by $H_{\cdot j}$ and $(i,j)$ element by $H_{ij}$ or $H_{i,j}$.

For real-valued functions $f_1$ and $f_2$ and $f_3$, we write 
\begin{equation*}
    f_1=\mathcal{O}(f_2)
\end{equation*}
if there exits a positive number $k_4$ such that $f_1\leq k_4f_2.$ We write
\begin{equation*}
    f_1 = \tilde{\mathcal{O}}_{f_3}(f_2)
\end{equation*}
if there exists a positive number $k_5$ such that $f_1\leq k_5f_2\times\text{poly}\log(f_3).$

\subsection{IPMs for LCQO} \label{sec: IPM}
In this work, LCQO is defined as follow.
\begin{Definition}[LCQO Problem]\label{def: LCQO}
For vectors $b\in \mathbb{R}^m$, $c\in \mathbb{R}^n$, and matrix $A\in \mathbb{R}^{m\times n}$ with $\rank(A)=m\leq n$, and symmetric positive semidefinite matrix $Q\in \mathbb{R}^{n \times n}$, we define the primal and dual LCQO problems as:
\begin{equation} \label{e:LCQP}
    (\text{P})\quad
    \begin{aligned}
    \min\  c^T&x + \frac{1}{2}x^TQx, \\
    {\rm s.t. }\;\;
    Ax &= b, \\
    x &\geq 0,
    \end{aligned}
    \qquad \qquad (\text{D}) \quad
    \begin{aligned}
    \max \  b^Ty - \frac{1}{2}x^TQ&x,\ \  \\
    {\rm s.t. }\;\;
    A^Ty + s - Qx&= c,\\
    s &\geq 0,
    \end{aligned}
\end{equation}
where $x\in \mathbb{R}^n$ is the vector of primal variables, and $y\in\mathbb{R}^m$, $s\in\mathbb{R}^n$ are vectors of the dual variables.
Problem $(P)$ is called the primal problem and $(D)$ is called the dual problem. 
\end{Definition}
The full-row-rankness of matrix $A$ implies that there is no all-zero row in matrix $A$. We further make the following assumption on matrix $A$.
\begin{Assumption}\label{assumption: no zero}
Matrix $A$ has no all-zero columns.
\end{Assumption}

\begin{Remark}
    When matrix $A$ has zero columns, without loss of generality, let us say the $n^{th}$ column is all-zero,  then we can introduce a new variable $x_{n+1}$ and rewrite the problem into
    \begin{equation*}
        \begin{aligned}
    \min\  \begin{bmatrix}
        c\\0
    \end{bmatrix}^T\begin{bmatrix}
        x\\x_{n+1}
    \end{bmatrix} + \frac{1}{2}\begin{bmatrix}
        x\\x_{n+1}
    \end{bmatrix}^T\begin{bmatrix}
        Q&0_{n\times1}\\0_{1\times n}&0
    \end{bmatrix}&\begin{bmatrix}
        x\\x_{n+1}
    \end{bmatrix}, \\
    {\rm s.t. }\;\;
    \begin{bmatrix}
        A_{\cdot 1} &\cdots & A_{\cdot (n-1)} & 0_{m\times 1} & 0_{m\times 1}\\
        0           &\cdots & 0               & 1     & -1
    \end{bmatrix}\begin{bmatrix}
        x\\
        x_{n+1}
    \end{bmatrix} &= \begin{bmatrix}
        b\\0
    \end{bmatrix}, \\
    x \geq 0,\ x_{n+1}&\geq 0.
    \end{aligned}
    \end{equation*}
    The new problem is equivalent to the original one. The new problem is still a LCQO problem and has fewer all-zero columns than the original problem. So we can repeat the procedure to eliminate all the all-zero columns.  In the worst case, we will get a new LCQO problem satisfying Assumption~\ref{assumption: no zero} with $2n-m$ variables and $n$ constraints.
\end{Remark}
\begin{Assumption}\label{assumption: ipc}
There exists a solution $(x,\ y,\ s)$ such that
\begin{equation*}
    Ax=b,\ x>0,\ A^Ty+s-Qx=c, \text{ and } s>0.
\end{equation*}
\end{Assumption}
The set of \textit{primal-dual feasible solutions} can be defined as
\begin{equation*}
    \mathcal{PD} \coloneqq \left\{(x,\ y,\ s)\in \mathbb{R}^n\times\mathbb{R}^m\times\mathbb{R}^n : Ax=b,\ A^Ty+s-Qx=c,\ (x,s)\geq0 \right\}
\end{equation*}
and the set of \textit{interior feasible primal-dual solutions} can be defined as
\begin{equation*}
    \mathcal{PD}^0 \coloneqq \left\{(x,\ y,\ s)\in \mathbb{R}^n\times\mathbb{R}^m\times\mathbb{R}^n : Ax=b,\ A^Ty+s-Qx=c,\ (x,s)>0 \right\}.
\end{equation*}
According to the strong duality, the set of optimal solutions can be defined as
\begin{equation*}
    \mathcal{PD}^\ast \coloneqq \left\{(x,\ y,\ s)\in \mathcal{PD} : xs=0 \right\},
\end{equation*}
where $xs$ denotes the Hadamard, i.e., componentwise product of $x$ and $s$. Let $\epsilon>0$, then the set of $\epsilon$-approximate solutions to Problem~\eqref{def: LCQO} can be defined as
\begin{equation}\label{def: approximate solution LCQO}
    \mathcal{PD}_{\epsilon} \coloneqq \left\{(x,\ y,\ s)\in \mathcal{PD} : x^Ts\leq n \epsilon \right\}.
\end{equation}
Let $X$ and $S$ be diagonal matrices of $x$ and $s$, respectively. Under Assumption~\ref{assumption: ipc}, for all $\mu>0$, the perturbed optimality conditions
\begin{equation}\label{syst: KKT}
    \begin{aligned}
    Ax &= b,\\
    A^T y + s - Qx &= c,\\
    XSe & = \mu e,\\
    (x,s) & \geq 0
    \end{aligned}
\end{equation}
have a unique solution $(x(\mu),\ y(\mu),\ s(\mu))$ that defines the primal and dual central path
\begin{equation*}
    \mathcal{CP} \coloneqq \left\{(x,y,s)\in \mathcal{PD}^0| x_is_i = \mu \text{ for } i\in\{1,\dots,n\}; \text{ for } \mu > 0  \right\}.
\end{equation*}
IPMs apply Newton's method to solve system (\ref{syst: KKT}). At each iteration of infeasible IPMs, a candidate solution to the primal-dual LCQO pair in \eqref{e:LCQP} is updated by solving the following linear system to find the Newton direction:
\begin{equation}\label{syst: full newton}
    \begin{bmatrix}
    A  & 0   & 0\\
    -Q & A^T & I\\
    S  & 0   & X
    \end{bmatrix} \begin{bmatrix}
    \Delta x\\\Delta y\\\Delta s
    \end{bmatrix}=\begin{bmatrix}
    r_p\\r_d\\r_c
    \end{bmatrix},
\end{equation}
where $(r_p,r_d,r_c)$ are residuals defined as
\begin{equation*}
    \begin{aligned}
    r_p &= b-Ax\\
    r_d &= c - A^Ty - s\\
    r_c &= \sigma\mu e - XSe,
    \end{aligned}
\end{equation*}
where $\sigma\in(0,1)$ is the barrier reduction parameter. If $r_p=0$ and $r_d=0$, then the solutions $(x,\ y,\ s)$ are primal and dual feasible. Alternatively, we can also define residuals in different ways as we will show later. Once the Newton direction is found, one can move along the direction but has to stay in a neighbourhood of the central path, which is defined at the end of this section.

When the linear system \eqref{syst: full newton} is solved inexactly, that actually leads to inexact infeasible IPMs. Many researchers have analyzed the performance of inexact infeasible IPMs (II-IPMs). For LCQO problems, \cite{lu2006iterative} propose an II-IPM using an iterative method to solve the Newton systems and obtain $\mathcal{O}(n^2\log(\frac{1}{\epsilon}))$ IPM iteration complexity. Here IPM iteration complexity does not include the complexity contributed by linear system solvers. However, it is known that feasible IPMs for LCQO problems can achieve $\mathcal{O}(\sqrt{n}\log(\frac{1}{\epsilon}))$ IPM iteration complexity \cite{kojima1989polynomial, monteiro1989interior, goldfarb1990n}. In \cite{gondzio2013convergence}, the author provides a general inexact feasible IPM for LCQO problems but has not discussed how to maintain feasibility when inexact linear system solvers are used. In this work, we will fill the gap by using a method inspired by some QIPM results \cite{mohammadisiahroudi2021inexact, augustino2021quantum} as we shall discuss later.

In this paper, we consider the following neighborhood of the central path
\begin{equation}\label{set:neighbourhood}
    \mathcal{N}_2 (\theta) \coloneqq \left\{(x,y,s)\in \mathcal{PD}^0|\|XSe-\mu e\|_2\leq\theta\mu  \right\},
\end{equation}
where $\theta \in (0, 1)$.

\subsection{Orthogonal Subspaces System}
Assuming that $(x,\ y,\ s) \in \mathcal{PD}^0$, to maintain the feasibility of the primal and dual variables, the first two linear equations in system~\eqref{syst: full newton} need to be solved with $r_p=0$ and $r_d=0$ exactly, which can be guaranteed if $\Delta x$ lies in the null space of $A$, denoted as $\Null(A)$, and $\Delta s= Q\Delta x - A^T \Delta y$. Accordingly, we can rewrite system~\eqref{syst: full newton} if we represent $\Delta x$ by a basis of $\Null(A)$. To do so, we can partition matrix $A$ to $A = \begin{bmatrix}
A_B&A_N
\end{bmatrix}$, where $A_B$ is a basis of $A$. Then we construct the following matrix
\begin{equation*}
    V = \begin{bmatrix}
    A_B^{-1}A_N \\ - I
    \end{bmatrix}.
\end{equation*}
Matrix $V$ has full column rank and satisfies $AV=0$, i.e., the columns of $V$ span the null space of $A$. Let $\Delta x = V\lambda$, where $\lambda\in \mathbb{R}^{n-m}$ is the unknown coefficient vector for $\Delta x$. Subsequently, we can rewrite system~\eqref{syst: full newton} by substituting $\Delta x$ and $\Delta s$ in the third equation as
\begin{equation}\label{syst: OSS}
    \begin{aligned}
    SV\lambda + X\left(QV\lambda - A^T \Delta y \right) = r_c \Leftrightarrow \begin{bmatrix}
    SV + XQV & - XA^T
    \end{bmatrix}\cdot \begin{bmatrix}
    \lambda\\\Delta y
    \end{bmatrix} = r_c.
    \end{aligned}
\end{equation}
A similar system was proposed and called "Orthogonal Subspaces System" (OSS) in \cite{augustino2021quantum, mohammadisiahroudi2021inexact} and we use the same name in this work. The matrix in the OSS system~\eqref{syst: OSS} is of size $n\times n$, and it is nonsingular. Even if the OSS system is solved inexactly, primal and dual feasibility is preserved by computing $\Delta x = V\lambda$ and $\Delta s= QV\lambda - A^T\Delta y$. Thus, we can conclude that residual will only show up in the third equation of~\eqref{syst: full newton}, i.e., $r_p=0$ and $r_d=0$. This nice property of the OSS system brings much convenience in the analysis of the proposed inexact IPM, and allows to prove the to-date best iteration complexity.  

\section{Inexact Feasible IPM with QLSAs}
In this section, we propose our IF-QIPM for LCQO problems. We start with the IF-IPM structure introduced by  \cite{gondzio2013convergence} and describe how to convert it into an IF-QIPM. Then we analyze the construction of the OSS system, and finally, we analyze the complexity for our IF-QIPM.

\subsection{IF-IPM for LCQO}
In \cite{gondzio2013convergence}, the author studies a general conceptual form IF-IPM for QCLO problems by assuming the feasibility of primal and dual variables, which induces the following system
\begin{equation}\label{syst: feasible full Newton}
    \begin{bmatrix}
    A  & 0   & 0\\
    -Q & A^T & I\\
    S  & 0   & X
    \end{bmatrix} \begin{bmatrix}
    \Delta x\\\Delta y\\\Delta s
    \end{bmatrix}=\begin{bmatrix}
    0\\0\\r_c
    \end{bmatrix},
\end{equation}
where $r_c = \sigma\mu e  - XSe$ with $\sigma\in(0,1)$ being the reduction factor of the central path parameter $\mu$, i.e., $\mu^{new} = \sigma\mu$. When system~\eqref{syst: feasible full Newton} is solved with $r_c = \sigma\mu e - XSe$ inexactly yielding an error $r$, if $\|r\|_2\leq \delta\|r_c\|_2$ for some $\delta\in(0,1)$, then the inexact IPM produces an $\epsilon$-approximate solution to Problem~\eqref{def: LCQO} in $\mathcal{O}(\sqrt{n}\log(1/\epsilon))$ iterations. The author of \cite{gondzio2013convergence} does not specify how to solve system~\eqref{syst: feasible full Newton} inexactly, how to preserve primal and dual feasibility, and how to   satisfy the convergence conditions described in \cite{gondzio2013convergence}. Specifically, the convergence conditions are posed on the right-hand-side and the inexactness error of the system (\ref{syst: feasible full Newton}).

Now we present a general procedure how to solve system~\eqref{syst: feasible full Newton} inexactly, while the inexactness error occurs only in the third equation of system~\eqref{syst: feasible full Newton}. Let $(\lambda, \Delta y)$ be an inexact solution for system~\eqref{syst: OSS} and $r$ be the corresponding inexactness error, so we have
\begin{equation*}
    \begin{aligned}
    \begin{bmatrix}
    SV + XQV & - XA^T
    \end{bmatrix}\cdot \begin{bmatrix}
    \lambda\\\Delta y
    \end{bmatrix} = r_c +r.
    \end{aligned}
\end{equation*}
The corresponding Newton step
\begin{equation*}
    \begin{aligned}
        \Delta x &= V\lambda\\
        \Delta s &= Q\Delta x - A^T \Delta y
    \end{aligned}
\end{equation*}
satisfies
\begin{equation*}
    \begin{bmatrix}
    A  & 0   & 0\\
    -Q & A^T & I\\
    S  & 0   & X
    \end{bmatrix}\cdot\begin{bmatrix}
    \Delta x\\\Delta y\\\Delta s
    \end{bmatrix}=\begin{bmatrix}
    0\\0\\r_c + r
    \end{bmatrix}.
\end{equation*}
Recall that once $(\lambda,\Delta y)$ is determined, then $(\Delta x, \Delta s)$ is also determined. An interesting property is that, if $(\lambda,\Delta y)$ and $(\Delta x, \Delta y, \Delta s)$ can be deduced from each other, then the OSS system and system~\eqref{syst: feasible full Newton} yield the same error term $r$. Hence the convergence conditions built upon system~\eqref{syst: feasible full Newton} can be directly examined using the residual $r_c$ and error $r$ of the OSS system. Let $\epsilon_{OSS}$ be the target accuracy of the OSS system~\eqref{syst: OSS}, i.e.,
\begin{equation*}
    \|(\lambda -  \lambda^\ast, \ \Delta  y - \Delta y^\ast\|_2\leq \epsilon_{OSS},
\end{equation*}
where $(\lambda^\ast, \ \Delta y^\ast)$ is the accurate solution. To make the IF-IPM converge, according to \cite{gondzio2013convergence}, we need
\begin{equation*}
    \begin{aligned}
        \|r\|_2  &= \left\| \begin{bmatrix}
                                SV + XQV & - XA^T
                                \end{bmatrix}\cdot \begin{bmatrix}
                                \lambda\\\Delta y
                                \end{bmatrix} - r_c\right\|_2\\
                   &\leq \left\|\begin{bmatrix}
                                SV + XQV & - XA^T
                                \end{bmatrix}\right\|_2 \epsilon_{OSS}\\
                    &\leq \delta \|r_c\|_2.
    \end{aligned}
\end{equation*}
So
\begin{equation*}
    \begin{aligned}
        \epsilon_{OSS}\leq \delta \frac{\|r_c\|_2}{\left\|\begin{bmatrix}
                                SV + XQV & - XA^T
                                \end{bmatrix}\right\|_2 }
    \end{aligned}
\end{equation*}
is sufficient for the IF-IPM to converge. We present the IF-IPM in Algorithm~\ref{alg: IF-IPM}. 
%
\begin{algorithm}[H] 
\caption{Short-step IF-IPM } 
\label{alg: IF-IPM}
\begin{algorithmic}[1]
 \STATE Choose $\epsilon >0$, $\delta\in(0,1)$, $\theta\in(0,1)$, $\beta\in(0,1)$ and $\sigma=(1-\frac{\beta}{\sqrt{n}})$. 
 \STATE $k \gets 0$
 \STATE Choose initial feasible interior solution $(x^0, y^0, s^0)\in \mathcal{N}(\theta)$ 
\WHILE {$(x^k,y^k,s^k)\notin \mathcal{PD}_\epsilon$}

\STATE $\mu^k \gets \frac{(x^k)^Ts^k}{n}$\label{alg-step:mu_k}
\STATE $\epsilon_{OSS}^k \gets \delta\|r_c^k\|_2/\left\|\begin{bmatrix}
                                S^kV + X^kQV^k & - X^kA^T
                                \end{bmatrix}\right\|_2$\label{alg-step:epsilon_k}
\STATE $(\lambda^k,\Delta y^k) \gets$ \textbf{solve} system~\eqref{syst: OSS} with accuracy $\epsilon_{OSS}^k$  \label{alg-step:LSA}
\STATE $\Delta x^k= V\lambda^k$ and $\Delta s^k=- A^T\Delta y^k$ \label{alg-step:x and s}
\STATE $(x^{k+1},y^{k+1},s^{k+1}) \gets (x^k,y^k,s^k)+(\Delta x^k,\Delta y^k,\Delta s^k)$\label{alg-step:update solution}
\STATE $k \gets k+1$
\ENDWHILE
\RETURN{$(x^k,y^k,s^k)$} 
\end{algorithmic}
\end{algorithm}
In the quantum-assisted IF-IPM, or IF-QIPM, we are proposing to accelerate Step \ref{alg-step:LSA} using quantum algorithms. In the next sections, we investigate how to use quantum algorithms to build and solve the OSS system and get the Newton direction.

\subsection{IF-QIPM for LCQO}
The pseudocode of our IF-QIPM is presented in Algorithm \ref{alg: IF-QIPM}. At each iteration of the IF-QIPM, we construct and solve system (\ref{syst: OSS}) and compute the Newton direction using quantum algorithms.
\begin{algorithm}[H] 
\caption{Short-step IF-QIPM } 
\label{alg: IF-QIPM}
\begin{algorithmic}[1]
 \STATE Choose $\epsilon >0$, $\delta\in(0,1)$, $\theta\in(0, \theta_0)$, $\beta\in(0,1)$ and $\sigma=(1-\frac{\beta}{\sqrt{n}})$. 
 \STATE $k \gets 0$
 \STATE Choose initial feasible interior solution $(x^0, y^0, s^0)\in \mathcal{N}(\theta)$ 
\WHILE {$(x^k,y^k,s^k)\notin \mathcal{PD}_\epsilon$}

\STATE $\mu^k \gets \frac{(x^k)^Ts^k}{n}$\label{alg-step:mu_k q}
\STATE $\epsilon_{OSS}^k \gets \delta\|r_c^k\|_2/\left\|\begin{bmatrix}
                                S^kV + X^kQV^k & - X^kA^T
                                \end{bmatrix}\right\|_2$\label{alg-step:epsilon_k q}
\STATE $(\lambda^k,\Delta y^k) \gets$ \textbf{solve} system~\eqref{syst: OSS} with accuracy $\epsilon_{OSS}^k$ quantumly\label{alg-step:QLSA}
\STATE $\Delta x^k= V\lambda^k$ and $\Delta s^k=- A^T\Delta y^k$ \label{alg-step:x and s q}
\STATE $(x^{k+1},y^{k+1},s^{k+1}) \gets (x^k,y^k,s^k)+(\Delta x^k,\Delta y^k,\Delta s^k)$\label{alg-step:update solution q}
\STATE $k \gets k+1$
\ENDWHILE
\RETURN{$(x^k,y^k,s^k)$}
\end{algorithmic}
\end{algorithm}
Here $\theta_0<1$ and its value will be discussed later. First, we introduce some notations to simplify the OSS system. In the $k^{th}$ iteration of Algorithm \ref{alg: IF-QIPM}, let
\begin{equation*}
    \begin{aligned}
    M^k = \begin{bmatrix}
    S^kV + X^kQV & - X^kA^T
    \end{bmatrix}, \ z^k = \begin{bmatrix}
    \lambda^k\\\Delta y^k
    \end{bmatrix}. 
    \end{aligned}
\end{equation*}
Then the OSS system can be rewritten as
\begin{equation*}
    M^kz^k=r_c^k.
\end{equation*}
As discussed in \cite{mohammadisiahroudi2021inexact}, to solve the OSS system~\eqref{syst: OSS} using quantum algorithms, we need to first rewrite it as the normalized Hermitian OSS system
\begin{equation}\label{syst: normalized hermitian OSS}
    \begin{aligned}
    \frac{1}{\sqrt{2}\left\|M^k\right\|_F}\begin{bmatrix}
    0   & M^k\\
    (M^k)^T & 0
    \end{bmatrix}\cdot \begin{bmatrix}
    0\\ z^k
    \end{bmatrix} = \frac{1}{\sqrt{2}\left\|M^k\right\|_F}.
    \end{aligned}\begin{bmatrix}
    r_c^k\\0
    \end{bmatrix}.
\end{equation}
To use the QLSAs mentioned earlier, we need to turn the linear system \eqref{syst: normalized hermitian OSS} into a quantum linear system using the block-encoding introduced in \cite{gilyen2018quantum}. To this end, we first decompose the coefficiennt matrix in linear system~\eqref{syst: normalized hermitian OSS} as
\begin{equation}\label{block encoding decomposition1}
    \begin{aligned}
    \frac{1}{\sqrt{2}\left\|M^k\right\|_F}\begin{bmatrix}
    0   & M^k\\
    (M^k)^T & 0
    \end{bmatrix} &= \frac{1}{\sqrt{2}\left\|M^k\right\|_F}\begin{bmatrix}
    0   & 0\\
    (M^k)^T & 0
    \end{bmatrix} + \frac{1}{\sqrt{2}\left\|M^k\right\|_F}\begin{bmatrix}
    0   & M^k\\
    0   & 0
    \end{bmatrix},
    \end{aligned}
\end{equation}
where
\begin{equation}\label{block encoding decomposition2}
    \begin{aligned}
    \begin{bmatrix}
    0   & 0\\
    (M^k)^T & 0
    \end{bmatrix} =&
    \begin{bmatrix}
        0_{n\times n} & 0_{n\times n} & 0_{n\times n}\\
        0_{(n-m)\times n} & V^T & 0_{(n-m)\times n}\\
        0_{m\times n} & 0_{m\times n} & -A
    \end{bmatrix}\times\\
    &\left(\begin{bmatrix}
        0_{n\times n} & 0_{n\times n}\\
        S^k & 0_{n\times n}\\
        0_{n\times n} & 0_{n\times n}
    \end{bmatrix} + 
    \begin{bmatrix}
        0_{n\times n} & 0_{n\times n} & 0_{n\times n}\\
        0_{n\times n} & Q & 0_{n\times n}\\
        0_{n\times n} & 0_{n\times n} & I_{n\times n}
    \end{bmatrix}
    \begin{bmatrix}
        0_{n\times n} & 0_{n\times n}\\
        X^k & 0_{n\times n}\\
        X^k & 0_{n\times n}
    \end{bmatrix}
    \right).
    \end{aligned}
\end{equation}
To compute matrix $V$, we need to find a basis matrix $A_B$ of matrix $A$ and we need to compute the inverse matrix $A_B^{-1}$. Both steps are nontrivial and can be expensive. However, we can reformulate the LCQO problem as follows
\begin{equation*}
    \begin{aligned}
    \min \ c^Tx &+ \frac{1}{2}x^TQx\\
    {\rm s.t. }  \ \begin{bmatrix}
        I &0 &A\\
        0 &I &-A
    \end{bmatrix}&\begin{bmatrix}
        x^\prime\\x^{\prime\prime}\\x
    \end{bmatrix} =\begin{bmatrix}
        b\\-b
    \end{bmatrix}\\
          x\geq0,\ x^\prime &\geq 0, x^{\prime\prime}\geq 0.
    \end{aligned}
\end{equation*}
In this case, we have an obvious basis 
\begin{equation*}
    A_B = \begin{bmatrix}
        I &0\\
        0 &I
    \end{bmatrix}
\end{equation*}
and matrix $V$ can be constructed efficiently
\begin{equation*}
    V = \begin{bmatrix}
        A_B^{-1} A_N\\-I
    \end{bmatrix}=\begin{bmatrix}
        \begin{bmatrix}
        I &0\\0 & I
    \end{bmatrix}\begin{bmatrix}
        A\\-A
    \end{bmatrix}\\-I
    \end{bmatrix} = \begin{bmatrix}
        A\\-A\\-I
    \end{bmatrix}.
\end{equation*}
Since matrix $A$ has no all-zero rows, matrix $V$ has no all-zero rows either. This property of the reformulation is useful in the analysis of the proposed IF-QIPM but we do not want to build the complexity analysis on the reformulated problem. So without loss of the generality we may make the following assumption. 
\begin{Assumption}\label{assumption: A}
Matrix $A$ is of the form
$A = \begin{bmatrix}
            I&A_N
        \end{bmatrix}$.
\end{Assumption}
To simplify the analysis, we further assume the input data are integers.
\begin{Assumption}\label{assumption: integer data}
    The input data of Problem~\eqref{def: LCQO} are integers.
\end{Assumption}
Following from the two assumptions above, we have the following lemma.
\begin{Lemma}\label{lemma: V}
Matrix $V$ equals to
\begin{equation*}
    V = \begin{bmatrix}
        A_N\\-I
    \end{bmatrix}
\end{equation*}
and
\begin{equation*}
    \min_{i=1,\dots, n}\{\|V_{i\cdot}\|_2^2\}  = \min\{1,\ \min_{i=1,\dots, m} \|(A_N)_{i\cdot}\|_2^2\}=1,
\end{equation*}
where $V_{i\cdot}$ and $(A_N)_{i\cdot}$ are the $i^{th}$ row of $V$ and $A_N$, respectively.
\end{Lemma}

%
%
Now we are ready to give $\theta_0$ in our definnition of central path neighbourhood, see~\eqref{set:neighbourhood}. We set
\begin{equation}\label{eq:theta}
    \theta_0 = \min \left\{\frac{1}{3\sqrt{n}},\ \frac{1}{4\left\|QVV^T\right\|_F + 1} \right\}.
\end{equation}
We also define $\omega^k$ as the maximum of the values of primal variables and dual slack variables in the $k^{th}$ iteration.
\begin{Definition}\label{def: omega}
Let $(x^k,y^k,s^k)$ be the a candidate solution for Problem~\eqref{e:LCQP}, then 
\begin{equation*}
    \omega^k = \max_{i\in \{1,\dots, n\}} \{x_i^k, s_i^k\}.
\end{equation*}
\end{Definition}
In this work, we assume access to quantum random access memory, QRAM. Then Step \ref{alg-step:QLSA} of Algorithm \ref{alg: IF-QIPM} consists of three parts: 1.) use block-encoding to build system~\eqref{syst: normalized hermitian OSS}; 2.) use QLSAs to solve system~\eqref{syst: normalized hermitian OSS}; 3.) use quantum tomography algorithms (QTAs) to extract classical solution. We use the block-encoding methods introduced in \cite{gilyen2018quantum} to block-encode linear system (\ref{syst: normalized hermitian OSS}). 
\begin{Proposition}
In the $k^{th}$ iteration of Algorithm \ref{alg: IF-QIPM}, use the block-encoding methods introduced in~\cite{gilyen2018quantum} and the decomposition described in equations (\ref{block encoding decomposition1}) and (\ref{block encoding decomposition2}), a
\begin{gather*}
    \Bigg(\frac{\sqrt{\|V\|_F^2 + \|A\|_F^2}\sqrt{2}\omega^k}{\|M^k\|_F}(\sqrt{2}\|Q\|_F+\sqrt{2}+1),\ \mathcal{O}(\poly\log(n)),\ \frac{\epsilon_{QLSA}}{\kappa_{M^k}^3}\Bigg)
\end{gather*}
-block-encoding of the matrix in the system (\ref{syst: normalized hermitian OSS}) can be implemented efficiently and the complexity will be dominated by the complexity of the QLSA step. Here $\epsilon_{QLSA}$ is the accuracy required for the QLSA step and $\kappa_{M^k}$ is the condition number of matrix $M^k$.
\end{Proposition}
\begin{proof}
See Appendix~\ref{sec: BE OSS} for proof.
\end{proof}
The complexity contributed by block-encoding is negligible compared with the complexity contributed by QLSAs and QTAs so we ignore it here. To establish the total complexity contributed by QLSAs and QTAs, we first need to analyze the accuracy of QLSA characterized by $\epsilon_{QLSA}$ and the accuracy of QTA characterized by $\epsilon_{QTA}$ and their relationship.

In each iteration, we use a QLSA to solve the block-encoded version of system~\eqref{syst: normalized hermitian OSS} and get an $\epsilon_{QLSA}$-approximate solution. Then we use a QTA to extract an $\epsilon_{QTA}$-approximate solution from the quantum machine. Here, for QLSA and QTA,  $\tilde{z}$ is an $\epsilon$-approximate solution of $z$ means
\begin{equation*}
    \left\| \frac{\Tilde{z} }{\|\Tilde{z} \|_2}- \frac{z}{\|z\|_2} \right\|_2 \leq \epsilon,
\end{equation*}
which is different from the concept of $\epsilon$-approximate solutions defined in~\eqref{def: approximate solution LCQO}.

Similar to \cite{mohammadisiahroudi2022efficient}, the QLSA we use is proposed by \cite{chakraborty2018power} and the QTA we use is proposed by \cite{van2022quantum}. Following the argument in Section 2 in \cite{mohammadisiahroudi2022efficient}, we can set the relationship among $\epsilon_{QLSA}$, $\epsilon_{QTA}$, and $\epsilon_{OSS}^k$ as
 \begin{equation}\label{accuracy: quantum}
    \begin{aligned}
    \epsilon_{QLSA} = \epsilon_{QTA} =\frac{1}{2}\cdot \frac{\sqrt{2}\|M^k\|_F}{\|r_c^k\|_2}\epsilon_{OSS}^k,
    \end{aligned}
\end{equation}
where $\epsilon_{OSS}^k$ is defined as the $\ell_2$ norm of the residual when solving system~\eqref{syst: normalized hermitian OSS} inexactly in the $k^{th}$ iteration. Here we did not add superscript for $\epsilon_{QLSA}$ and $\epsilon_{QTA}$ and the reason shall be revealed later. Let
\begin{equation*}
    \begin{bmatrix}
        \tilde{0}^k\\
        \tilde{z}^k
    \end{bmatrix}
\end{equation*}
be an inexact solution for system~\eqref{syst: normalized hermitian OSS} in the $k^{th}$ iteration. Then the norm of residual of system~\eqref{syst: normalized hermitian OSS}, which is $\epsilon_{OSS}^k$, and the norm of residual of system~\eqref{syst: OSS}, which is $\|M^k\tilde{z}^k-r_c^k\|_2$, satisfies
\begin{align*}
    \epsilon_{OSS}^k &= \left\| \frac{1}{\sqrt{2}\|M^k\|_F}\begin{bmatrix}
        0 & M^k\\ (M^k)^T & 0
    \end{bmatrix}\begin{bmatrix}
        \tilde{0}^k\\\tilde{z}^k 
    \end{bmatrix} - \frac{1}{\sqrt{2}\|M^k\|_F}\begin{bmatrix}
        r_c^k\\0
    \end{bmatrix} \right\|_2\\
    &= \left\| \frac{1}{\sqrt{2}\|M^k\|_F}\begin{bmatrix}
        M^k\tilde{z}^k\\(M^k)^T\tilde{0}^k 
    \end{bmatrix} - \frac{1}{\sqrt{2}\|M^k\|_F}\begin{bmatrix}
        r_c^k\\0
    \end{bmatrix} \right\|_2\\
    &\geq \left\| \frac{1}{\sqrt{2}\|M^k\|_F}M^k\tilde{z}^k - \frac{1}{\sqrt{2}\|M^k\|_F}r_c^k \right\|_2\\
    &\geq \frac{1}{\sqrt{2}\|M^k\|_F} \|M^k\tilde{z}^k-r_c^k\|_2.
\end{align*}
Recall that the error arising from the OSS system~\eqref{syst: OSS} is the same as the error in the full Newton system~\eqref{syst: feasible full Newton}, then we can directly use the convergence condition provided in Gondzio's analysis to the IF-IPM scheme in \cite{gondzio2013convergence}, i.e.,
\begin{equation*}
    \|M^k\tilde{z}^k-r_c^k\|_2\leq\delta \|r_c^k\|_2,
\end{equation*}
where $\delta\in(0,1)$ is a constant parameter. We can require
$$
\|M^k\tilde{z}^k-r_c^k\|_2\leq \sqrt{2}\|M^k\|_F \epsilon_{OSS}^k\leq \delta \|r_c^k\|_2
$$
and it follows that
\begin{equation*}
    \epsilon_{OSS}^k\leq \frac{\delta\|r_c^k\|_2}{\sqrt{2}\left\|M^k\right\|_F}.
\end{equation*}
Then choosing
\begin{equation*}
    \epsilon_{QLSA} = \epsilon_{QTA} = \frac{\|M^k\|_F\epsilon_{OSS}^k}{\sqrt{2}\|r_c^k\|_2} \leq \frac{\delta}{2}
\end{equation*}
ensures the convergence of the IF-QIPM. The complexities for each step are also available now. Using the QLSA from \cite{chakraborty2018power} and QTA from \cite{van2022quantum}, we have the complexity for QLSA and QTA
\begin{equation*}
    \begin{aligned}
    T_{QLSA} &= \tilde{\mathcal{O}}_{n,\Bar{\omega},\frac{1}{\epsilon}}\left(\kappa_{M^k} \frac{\omega^k}{\|M^k\|_F} \right),\\
    T_{QTA}  &= \tilde{\mathcal{O}}_{n}\left(n\right).
    \end{aligned}
\end{equation*}
Note that the complexity of the block-encoding procedure is dominated by that of QLSA and QTA and thus we ignore the complexity contributed by  block-encoding. In Step \ref{alg-step:x and s q}, the complexity contributed by computing Newton step from OSS solution is $\mathcal{O}(n^2)$. The total complexity for the $k^{th}$ iteration of IF-QIPM  will be
\begin{equation}\label{complexity}    
\tilde{\mathcal{O}}_{n,\Bar{\omega},\frac{1}{\epsilon}}\left(\frac{n\omega^k\kappa_{M^k}}{\|M^k\|_F} + n^2\right).
\end{equation}

\subsubsection{Bound for $\omega^k/\|M^k\|_F$}\label{sec: bound M}
In this section, all the quantities we consider are from the $k^{th}$ iteration. For simplicity,  we ignore superscript $k$ in this section unless we need it. Using the property of trace, we have
\begin{equation*}
\begin{aligned}
    \|M\|_F^2 &= \trace(M^TM)\\
              &= \trace\left((SV + XQV)(SV + XQV)^T + XA^TAX  \right)\\
              &= \trace\left((SV + XQV)(SV + XQV)^T\right) + \trace\left(XA^TAX  \right)\\
              &= \trace\left(SVV^TS\right) + \trace\left(XQVV^TS\right) + \trace\left(SVV^TQX\right) + \trace\left(XQVV^TQX\right) + \trace\left(XA^TAX  \right).
\end{aligned}
\end{equation*}
For the non-symmetric term, due to cyclic invariant property of trace, we have
\begin{equation*}
    \trace\left(XQVV^TS\right) = \trace\left(SXQVV^T\right).
\end{equation*}
Recall the central path neighbourhood we defined in (\ref{set:neighbourhood}), we define a matrix $E$ such that
\begin{equation}\label{matrix E}
    E = \frac{1}{\mu\theta} (XS - \mu I). 
\end{equation}
It is obvious that $E$ is a diagonal matrix and satisfies
\begin{equation*}
    \|Ee\|_2 < 1,
\end{equation*}
which leads to
\begin{equation*}\label{matrix E property}
    |\trace(E)|\leq \|Ee\|_1\leq \sqrt{n}\|E\|_F = \sqrt{n}\|Ee\|_2 <\sqrt{n} \text{ and } I-E\succ 0 \text{ and } I+E\succ 0.
\end{equation*}
With this, we can have
\begin{equation*}
\begin{aligned}
    \trace\left(XQVV^TS\right) &= \trace\left(SXQVV^T\right)\\
                               &= \trace\left((\theta \mu E + \mu I)QVV^T\right)\\
                               &= \trace\left(\theta \mu E QVV^T\right) + \trace\left(\mu QVV^T \right).
\end{aligned}
\end{equation*}
For the second term, we know $Q$ and $V^TQV$ are both positive semidefinite. So we can have
\begin{equation*}
    \trace\left(QVV^T\right) = \trace\left(V^TQV\right)\geq 0
\end{equation*}
because of the cyclic invariant property of trace. According to the Cauchy–Schwarz inequality, we have
\begin{equation*}
    \trace\left(E QVV^T\right)^2 \leq \|E\|_F^2 \|QVV^T\|_F^2.
\end{equation*}
So we have
\begin{equation*}
    \trace\left(EQVV^T\right) \geq -\|QVV^T\|_F.
\end{equation*}
Thus, we have
\begin{equation*}
\begin{aligned}
    \trace\left(XQVV^TS\right) &= \trace\left(\theta \mu E QVV^T\right) + \trace\left(\mu QVV^T \right)\\
    &\geq \mu\left(\trace\left(QVV^T\right) - \theta\|QVV^T\|_F \right)\\
    &\geq -\theta\mu\|QVV^T\|_F\\
    &\geq - \frac{\mu}{4},
\end{aligned}
\end{equation*}
where the last inequality holds due to condition~\eqref{eq:theta}. So we can bound $\|M\|_F$ by
\begin{equation*}
    \begin{aligned}
    \|M\|_F^2 &= \trace\left(SVV^TS\right) + \trace\left(XQVV^TS\right) + \trace\left(SVV^TQX\right) + \trace\left(XQVV^TQX\right) + \trace\left(XA^TAX  \right)\\
              &\geq \trace\left(SVV^TS\right) + \trace\left(XQVV^TQX\right) + \trace\left(XA^TAX  \right) - \frac{\mu}{2}.
    \end{aligned}
\end{equation*}
Since $XQVV^TQX\succeq 0$, we have
\begin{equation*}
    \begin{aligned}
    \|M\|_F^2 &\geq \trace\left(SVV^TS\right) + \trace\left(XA^TAX  \right) - \frac{\mu}{2}.
    \end{aligned}
\end{equation*}
Since $X$ and $S$ are both positive diagonal matrices, we have
\begin{equation*}
    \begin{aligned}
    \|M\|_F^2 &\geq \trace\left(SVV^TS\right) + \trace\left(XA^TAX  \right) - \frac{\mu}{2}\\
    &= \sum_{i}s_i^2(VV^T)_{ii} + \sum_{i}x_i^2 (A^TA)_{ii} -\frac{\mu}{2}\\
    &\geq \omega^2 - \frac{\mu}{2}.
    \end{aligned}
\end{equation*}
As we said in the very beginning of this section, at each iteration $\omega$ is indeed $\omega^k$ but the superscript is ignored here. Now we are going to find a bound for $\mu$ so we can further bound $\|M\|_F^2$. Since $\omega$ is the upper bound for the magnitude of the primal and dual slack variables, we have
\begin{align*}
    \omega^2 &\geq x_is_i.
\end{align*} 
Recall the definition of matrix $E$, see \eqref{matrix E}. So we have
\begin{equation*}
    \begin{aligned}
    \omega^2 &\geq x_is_i = \mu + \theta \mu E_{ii} \geq \mu - \theta \mu = (1-\theta)\mu.
\end{aligned}
\end{equation*}
So
\begin{equation*}
    \begin{aligned}
    \|M\|_F^2 &\geq \omega^2  - \frac{\mu }{2} \geq \omega^2  - \frac{1}{2}\frac{\omega^2}{1-\theta} \geq \omega^2  - \frac{1}{2}\frac{\omega^2}{1-1/3} =\frac{\omega^2 }{4},
    \end{aligned}
\end{equation*}
where the last inequality follows from the bound for $\theta$, see~\eqref{eq:theta}. So we have
\begin{equation*}
\begin{aligned}
    \frac{\omega}{\|M\|_F} \leq 2=\mathcal{O}\left(1\right).
\end{aligned}
\end{equation*}

\subsubsection{Bound for $\kappa_{M^k}$}\label{sec: kappa}
Similar to the previous section, we ignore the supercript $k$ unless we need it. We will start with a general result and then work on the matrix $M^k$. The following lemma is a well-known result regarding condition numbers of matrices and can be proven using Courant-Fischer-Weyl Min-Max principle \cite{horn2012matrix}.
\begin{Lemma}\label{lemma:singular value QVQ}
For any full row rank matrix $P\in \mathbb{R}^{m\times n}$ and symmetric positive definite matrix $D\in \mathbb{R}^{n\times n}$, their condition number satisfies
$$
\kappa(PD P^T) \leq \kappa(D)\kappa(PP^T).
$$
\end{Lemma}

Next, we analyze the matrix in the OSS system~\eqref{syst: normalized hermitian OSS}. Specifically, we focus on $M^TM$ since we are interested in the spectral property of the OSS system~\eqref{syst: normalized hermitian OSS}. Using the matrix $E$ defined in~\eqref{matrix E},  we have the following decomposition
\begin{equation*}
\begin{aligned}
M^TM &= \begin{bmatrix}
V^T(S+XQ)^T(S+XQ)V&-V^T(S+XQ)^TXA^T\\-AX(S+XQ)V& AX^2A^T
\end{bmatrix}\\
&= \begin{bmatrix}
V^T(S+XQ)^T(S+XQ)V&-V^T\mu\left(\theta E\right)A^T - V^TQ^TX^2A^T\\-A\mu\left(\theta E\right)V - AX^2QV& AX^2A^T
\end{bmatrix}\\
&= \begin{bmatrix}
V^T&0\\0&A
\end{bmatrix}\begin{bmatrix}
(S+XQ)^T(S+XQ)&-\mu\theta E-QX^2\\-\mu\theta E-X^2Q&X^2
\end{bmatrix}\begin{bmatrix}
V^T&0\\0&A
\end{bmatrix}^T.
\end{aligned}
\end{equation*}
The second equality holds because
\begin{equation*}
\begin{aligned}
    -V^TSXA^T - V^TQ^TX^2A^T &= -V^T\mu\left(I + \theta E\right)A^T - V^TQ^TX^2A^T\\
    &= -V^T\mu\left(\theta E\right)A^T - V^TQX^2A^T.
\end{aligned}
\end{equation*}
Here we used that $AV=0$ and $Q$ is symmetric. Then, plugging~\eqref{matrix E} into the first diagonal block of the decomposition we obtained earlier, we have
\begin{adjustwidth}{-\extralength}{0cm}
\begin{equation}
\begin{aligned}
M^TM &= \begin{bmatrix}
V^T&0\\0&A
\end{bmatrix}\left(\begin{bmatrix}
S^2 + 2\mu Q +\mu\theta (EQ+QE)+ QX^2Q & -\mu\theta E-QX^2\\-\mu\theta E-X^2Q&X^2
\end{bmatrix}\right)\begin{bmatrix}
V^T&0\\0&A
\end{bmatrix}^T\\
&= \begin{bmatrix}
V^T&0\\0&A
\end{bmatrix}\left(\begin{bmatrix}
    S^2 + 2\mu Q + \mu\theta (EQ+QE) & -\mu\theta E\\-\mu\theta E&0
\end{bmatrix} + \begin{bmatrix}
QX^2Q & -QX^2\\-X^2Q&X^2
\end{bmatrix}\right)\begin{bmatrix}
V^T&0\\0&A
\end{bmatrix}^T\\
&= \begin{bmatrix}
V^T&0\\0&A
\end{bmatrix}\left(\begin{bmatrix}
    I & -Q\\ 0 & I
\end{bmatrix}\begin{bmatrix}
    S^2 + 2\mu Q  & -\mu\theta E\\-\mu\theta E&0
\end{bmatrix}\begin{bmatrix}
    I &0\\-Q &I
\end{bmatrix} + \begin{bmatrix}
    I & -Q\\ 0 & I
\end{bmatrix}\begin{bmatrix}
0 & 0\\0 &X^2
\end{bmatrix}\begin{bmatrix}
    I &0\\-Q &I
\end{bmatrix}\right)\begin{bmatrix}
V^T&0\\0&A
\end{bmatrix}^T\\
&= \begin{bmatrix}
V^T&0\\0&A
\end{bmatrix}\begin{bmatrix}
    I & -Q\\ 0 & I
\end{bmatrix}\begin{bmatrix}
    S^2 + 2\mu Q  & -\mu\theta E\\-\mu\theta E&X^2
\end{bmatrix}\begin{bmatrix}
    I &0\\-Q &I
\end{bmatrix}\begin{bmatrix}
V^T&0\\0&A
\end{bmatrix}^T.
\end{aligned}
\end{equation}
\end{adjustwidth}
The first two matrices are nonsingular, so we can apply the Lemma \ref{lemma:singular value QVQ} and thus we only need to study the middle matrix. Denote the middle matrix by $\Psi$. Observe that $\Psi$ is almost the same as its counterpart in \cite{mohammadisiahroudi2021inexact}. Subsequently we have the following result regarding the spectral property of $M^k$.
\begin{Lemma}\label{lemma: kappa of M}
When $(x,y,s)\in \mathcal{N}(\theta)$ and $\theta\in \left(0,\ \min \left\{\frac{1}{3\sqrt{n}},\ \frac{1}{4\|QVV^T\|_F + 1}\right\} \right)$, the condition number of matrix $M^k$ satisfies
$$
\kappa_{M^k} = \mathcal{O}\left(\frac{(\omega^k)^2 +  \mu^k\sigma_{\max}(Q)}{\mu^k}\kappa_{VAQ}  \right),
$$
where $\kappa_{VAQ}$ is the condition number of the matrix $\begin{bmatrix}
    V^T&0\\0&A
\end{bmatrix}\begin{bmatrix}
    I &-Q\\ 0&I
\end{bmatrix}$.
\end{Lemma}
\begin{proof}
The proof is in Appendix \ref{sec: spectral analysis for the middle matrix}.
\end{proof}
Putting all these together, we have the complexity for our IF-QIPM for LCQO problems.
\begin{Theorem}\label{theorem: main}
The IF-QIPM for LCQO problems stops with final duality gap less than $\epsilon$ in at most $\mathcal{O}\left( \sqrt{n}\log(1/\epsilon)\right)$ IPM iterations and in each IPM iteration, the Newton direction can be obtained with complexity  $\tilde{\mathcal{O}}_{n,\Bar{\omega},\frac{1}{\epsilon}}\left(n\left(\frac{\bar{\omega}^2}{\epsilon}+ \sigma_{\max}(Q)\right)\kappa_{VAQ} + n^2\right)$, where $\Bar{\omega} = \max_{k}\omega^k$.
\end{Theorem}
\begin{proof}
The complexity bound for the IPM iterations comes from the result in \cite{gondzio2013convergence}. According to \eqref{complexity}, the complexity for obtaining the Newton direction is
\begin{align}
    \tilde{\mathcal{O}}_{n,\Bar{\omega},\frac{1}{\epsilon}}\left(\frac{n\omega^k\kappa_{M^k}}{\|M^k\|_F} + n^2\right).
\end{align}
Combining this with the result in Sec.~\ref{sec: bound M}, the bound in Lemma \ref{lemma: kappa of M},  and $\mu^k\geq \epsilon$, we have
\begin{align}\tilde{\mathcal{O}}_{n,\Bar{\omega},\frac{1}{\epsilon}}\left(\frac{n\omega^k\kappa_{M^k}}{\|M^k\|_F} + n^2\right) = \tilde{\mathcal{O}}_{n,\Bar{\omega},\frac{1}{\epsilon}}\left(n\left(\frac{\bar{\omega}^2}{\epsilon}+ \sigma_{\max}(Q)\right)\kappa_{VAQ} + n^2\right).
\end{align}
\end{proof}

\section{Application in Support Vector Machine Problems}
In this section, we discuss how to use our IF-QIPM to solve SVM problems. We show that our algorithm can solve $l_1$-norm soft margin SVM problems with best complexity compared with any existing classical or quantum algorithms.

The ordinary SVM problem works on a linearly separable dataset, in which the data points have binary labels. The ordinary SVM aims to find a hyperplane correctly separating the data points with maximum margin. However, in practice the data points are not necessarily linearly separable. To allow mislabelling, the concept of soft margin SVM was introduced in \cite{cortes1995support}. Let $\{(\phi_i,\zeta_i)\in\mathbb{R}^m\times\{-1,+1\}|i=1,\dots,n\}$ be the set of data points, $\Phi$ be a matrix with $i$th column being $\phi_i$, and $Z$ be a diagonal matrix with $i$th diagonal element being $\zeta_i$. The SVM problem with $l_1$-norm soft margin can be formulated as below.
\begin{equation}\label{SVM: l1}
    \begin{aligned}
    \min_{(\xi,w,t)\in\mathbb{R}^n\times\mathbb{R}^m\times\mathbb{R}}\ &\frac{1}{2}\|w\|_2^2 + C \|\xi\|_1\\
              &\zeta_i \left(\langle w, \phi_i\rangle + t\right) \geq 1-\xi_i, \ i=1,\dots,n\\
                    &\xi_i\geq 0,\ i= 1,\dots,n.
    \end{aligned}
\end{equation}
Here $(w,t)$ determines a hyperplane and $C$ is a penalty parameter. In \cite{kerenidis2021quantum}, the authors rewrote the SVM problem as a second order conic optimization (SOCO) problem and use the quantum algorithm they proposed to solve the resulting SOCO problem. They claim the complexity of their algorithm has $\mathcal{O}(n^2)$ dependence on the dimension, which is better than any classical algorithm. However, the algorithm in \cite{kerenidis2021quantum} is invalid. Their algorithm is an Inexact Infeasible-QIPM (II-QIPM) while they used the IPM complexity for Feasible-QIPM, which ignores at least $\mathcal{O}(n^{1.5})$ dependence on $n$. They also missed the symmetrization of the Newton step, which is necessary for SOCO problems and makes their Newton step invalid. 

Aside from \cite{kerenidis2021quantum}, some pure quantum algorithms for SVM problems are also proposed. In \cite{rebentrost2014quantum}, the authors propose a pure quantum algorithm for SVM problems. They claim the complexity is $\mathcal{O}(\kappa_{\text{eff}}^3 \epsilon^{-3}\log(mn))$, where $\kappa_{\text{eff}}$ is the condition number of a matrix involving the kernel matrix and $\epsilon$ is the accuracy. In the worst case, $\kappa_{\text{eff}}=\mathcal{O}(m)$. Their complexity is worse than ours regarding the dependence of dimension and accuracy. In addition, their algorithm does not provide classical solutions. Namely, the solution is in the quantum machine and we can not read or use it in a classical computer. However, our algorithm produces a classical solution.

To convert the problem into standard form LCQO, we introduce $(w_+, w_-)\in\mathbb{R}_+^m\times\mathbb{R}_+^m$, $(t_+, t_-)\in\mathbb{R}_+\times\mathbb{R}_+$, and a slack variable $\rho\in\mathbb{R}_+^n$. Then we can get the following formulation
\begin{equation*}
    \begin{aligned}
    \min_{w_+,w_-,t_+,t_-,\xi,\rho} \ &\frac{1}{2}\|w_+-w_-\|_2^2 + C \|\xi\|_1\\
              &\zeta_i\left(\left<w_+-w_-,\phi_i\right>+t_+-t_-\right) + \xi_i -\rho_i =1, \ i=1,\dots,n\\
                    &(\xi,w_+,w_-,t_+,t_-,\rho)\geq0.
    \end{aligned}
\end{equation*}
It is a standard form LCQO problem with nonnegative variables $(w_+,w_-,t_+,t_-,\xi,\rho)\in \mathbb{R}^m\times\mathbb{R}^m\times\mathbb{R}\times\mathbb{R}\times\mathbb{R}^n\times\mathbb{R}^n$ and parameters
\begin{equation*}
    \begin{aligned}
    c &= \begin{bmatrix} 0_{2m+2}\\Ce_{n}\\ 0_{n}  \end{bmatrix}\\
    Q &= \begin{bmatrix}
            I_{m\times m} & -I_{m\times m} & 0_{m\times(2+2n)}\\
            -I_{m\times m} & I_{m\times m} & 0_{m\times(2+2n)}\\
            0_{(2+2n)\times m} & 0_{(2+2n)\times m} & 0_{(2+2n)\times(2+2n)}
         \end{bmatrix}\\
    A &= \begin{bmatrix}Z \Phi^T & - Z \Phi^T & Z & -Z & I_{n\times n} &-I_{n\times n} \end{bmatrix}\\
    b &= e.
    \end{aligned}
\end{equation*}
So we can use the proposed IF-QIPM for LCQO problems to solve the $\ell_1$-norm soft margin SVM problems and get an $\epsilon$-approximate solution with complexity $$\tilde{\mathcal{O}}_{n,\Bar{\omega},\frac{1}{\epsilon}}\left(n^{1.5}\left(\frac{\bar{\omega}^2}{\epsilon}+ \sigma_{\max}(Q)\right)\kappa_{VAQ} + n^{2.5}\right).$$
This dependence on dimension is better than any existing quantum or classical algorithm.

\section{Discussion}
In this work, we present an IF-QIPM for LCQO problems by combining the IF-IPM framework proposed in \cite{gondzio2013convergence} and the OSS system introduced in \cite{mohammadisiahroudi2021inexact}. Our algorithm has $n^{1.5}$ dependence on $n$, which is better than any existing algorithms for LCQO problems. The dependence on the accuracy is polynomial, which is worse than classic IPMs. Iterative refinement method might help improve the dependence on the accuracy but that could be another work. 

\vspace{6pt} 



\authorcontributions{Conceptualization, Zeguan Wu and Tam\'as Terlaky; Methodology, Zeguan Wu; Supervision, Xiu Yang and Tam\'as Terlaky; Validation, Zeguan Wu, Mohammadhossein Mohammadisiahroudi, Brandon Augustino, Xiu Yang and Tam\'as Terlaky; Writing – original draft, Zeguan Wu; Writing – review \& editing, Zeguan Wu, Mohammadhossein Mohammadisiahroudi, Brandon Augustino, Xiu Yang and Tam\'as Terlaky.
}


\funding{This work was supported by Defense Advanced Research Projects Agency as part of the project W911NF2010022: {\em The Quantum
Computing Revolution and Optimization: Challenges and Opportunities}.}

\institutionalreview{
Not applicable
.}

\informedconsent{Not applicable
%
.}

\dataavailability{Not applicable.
} 


\conflictsofinterest{The funder had no role in the design of the study; in the writing of the manuscript; or in the decision to publish the~results.}



\abbreviations{Abbreviations}{
The following abbreviations are used in this manuscript:\\

\noindent 
\begin{tabular}{@{}ll}
IF-IPM  &Inexact Feasible Interior Point Method\\
IF-QIPM & INexact Feasible Quantum Interior Point Methods\\
IPM & Interior Point Method\\
LCQO    &Linearly Constrained Quadratic Optimization\\
LO  & Linear Optimization\\
OSS & Orthogonal Subspace System\\
QIPM & Quantum Interior Point Method\\
QLSA & Quantum Linear System Algorithm\\
QTA & Quantum Tomography Algorithm\\
SOCO    & Second Order Conic Optimization\\
SVM &Support Vector Machine
\end{tabular}
}

\appendixtitles{yes} 
\appendixstart
\appendix
\section[\appendixname~\thesection]{Block-encoding of the OSS system}\label{sec: BE OSS}
%
In this section, we ignore the superscript $k$ for simplicity. As described in  Eq.~\eqref{block encoding decomposition1}, we first block encode each of the matrices involved in (\ref{block encoding decomposition2}). With $V,\ A,\ S$ and $X$ given and are stored in a quantum accessible data structure (we ignore the complexity to store the classical information into the quantum machine). For the first matrix
\begin{equation*}
    M_1=\begin{bmatrix}
        0_{n\times n} & 0_{n\times n} & 0_{n\times n}\\
        0_{(n-m)\times n} & V^T & 0_{(n-m)\times n}\\
        0_{m\times n} & 0_{m\times n} & -A
    \end{bmatrix},
\end{equation*}
a $$\left(\sqrt{\|V\|_F^2 + \|A\|_F^2}, \mathcal{O}(\poly\log(n)), \epsilon_1 \right)$$-block-encoding of $M_1$ can be implemented according to Lemma 50 from \cite{gilyen2018quantum} efficiently. 

The second matrix
\begin{equation*}
    M_2 = \begin{bmatrix}
        0_{n\times n} & 0_{n\times n}\\
        S & 0_{n\times n}\\
        0_{n\times n} & 0_{n\times n}
    \end{bmatrix}
\end{equation*}
is both 1-row-sparse and 1-column-sparse. By the definition of $\omega$, each element of $M_2/\omega$ has absolute value at most 1. According to Lemma 48 in \cite{gilyen2018quantum}, a $$\left(1, \mathcal{O}(\poly\log(n)), \epsilon_2 \right)$$-block-encoding of $M_2/\omega$ can be implemented efficiently.

The third matrix
\begin{equation*}
    M_3 = \begin{bmatrix}
        0_{n\times n} & 0_{n\times n} & 0_{n\times n}\\
        0_{n\times n} & Q & 0_{n\times n}\\
        0_{n\times n} & 0_{n\times n} & I_{n\times n}
    \end{bmatrix}
\end{equation*} can be decomposed into
$$
M_3 = \begin{bmatrix}
        0_{n\times n} & 0_{n\times n} & 0_{n\times n}\\
        0_{n\times n} & Q & 0_{n\times n}\\
        0_{n\times n} & 0_{n\times n} & 0_{n\times n}
    \end{bmatrix} + \begin{bmatrix}
        0_{n\times n} & 0_{n\times n} & 0_{n\times n}\\
        0_{n\times n} & 0_{n\times n} & 0_{n\times n}\\
        0_{n\times n} & 0_{n\times n} & I_{n\times n}
    \end{bmatrix}.
$$
Then we can block-encode the two matrices first, and  then apply linear combination to obtain $M_3$. In fact, a $$\left(\|Q\|_F, \mathcal{O}(\poly\log(n)), \epsilon_3\right)$$-block-encoding of the left matrix can be implemented according to Lemma 50 from \cite{gilyen2018quantum} efficiently
and a 
\begin{equation*}
    \begin{aligned}
    \left(1, \mathcal{O}(\poly\log(n)), \epsilon_3 \right)
\end{aligned}
\end{equation*}
-block-encoding of the right matrix can be implemented efficiently according to Lemma 48 in \cite{gilyen2018quantum}. With the state-preparation cost of the linear combination coefficient vector $(1, 1)$ neglected, a
$$
\left(\|Q\|_F+1, \mathcal{O}(\poly\log(n)), (\|Q\|_F+1)\epsilon_3\right)
$$
-block-encoding of $M_3$ can be implemented efficiently according to Lemma 52 from \cite{gilyen2018quantum}.

The fourth matrix
\begin{equation*}
    M_4 = \begin{bmatrix}
        0_{n\times n} & 0_{n\times n}\\
        X & 0_{n\times n}\\
        X & 0_{n\times n}
    \end{bmatrix}
\end{equation*}
is 1-row-sparse and 2-column-sparse. After being scaled by $\frac{1}{\omega}$, each element of $M_4/\omega$ has absolute value at most 1. According to Lemma 48 in \cite{gilyen2018quantum}, a $$\left(\sqrt{2}, \mathcal{O}(\poly\log(n)), \epsilon_4 \right)$$-block-encoding of $M_4/\omega$ can be implemented efficiently. \\

For the matrix multiplication $M_3M_4/\omega$, a $$\left(\sqrt{2}\|Q\|_F+\sqrt{2}, \mathcal{O}(\poly\log(n)), (\|Q\|_F+1)(\sqrt{2}\epsilon_3 + \epsilon_4) \right)$$-block-encoding can be implemented efficiently according to Lemma 53 from \cite{gilyen2018quantum}. \\

For the linear combination $M_2/\omega + M_3M_4/\omega$, the cost for the state-preparation of the coefficient vector $(1,1)$ is negligible and thus a $$\left(\sqrt{2}\|Q\|_F+\sqrt{2} +1, \mathcal{O}(\poly\log(n)), (\sqrt{2}\|Q\|_F+\sqrt{2} +1)(\epsilon_3 + \frac{1}{\sqrt{2}}\epsilon_4)  \right)$$-block-encoding can be implemented efficiently according to Lemma 52 from \cite{gilyen2018quantum}. 

For the matrix multiplication of $M_1(M_2/\omega + M_3M_4/\omega)$, a
\begin{gather*}
    \Bigg(\sqrt{\|V\|_F^2 + \|A\|_F^2}(\sqrt{2}\|Q\|_F+\sqrt{2}+1),\\ \mathcal{O}(\poly\log(n)),\\ \sqrt{\|V\|_F^2 + \|A\|_F^2} (\sqrt{2}\|Q\|_F+\sqrt{2}+1)(\epsilon_3 + \frac{1}{\sqrt{2}}\epsilon_4)  + (\sqrt{2}\|Q\|_F+\sqrt{2}+1)\epsilon_1 \Bigg)
\end{gather*}
-block-encoding can be implemented efficiently according to Lemma 53 from \cite{gilyen2018quantum}.

Finally, considering that the complexity of state-preparation of the vector $$(\frac{\omega}{\sqrt{2}\|M\|_F}, \frac{\omega}{\sqrt{2}\|M\|_F})$$ can be neglected, a
\begin{gather*}
    \Bigg(\frac{\sqrt{\|V\|_F^2 + \|A\|_F^2}\sqrt{2}\omega}{\|M\|_F}(\sqrt{2}\|Q\|_F+\sqrt{2}+1),\\ \mathcal{O}(\poly\log(n)),\\ \frac{\sqrt{\|V\|_F^2 + \|A\|_F^2}\sqrt{2}\omega}{\|M\|_F}(\sqrt{2}\|Q\|_F+\sqrt{2}+1)^2\left(\sqrt{\|V\|_F^2 + \|A\|_F^2}(\epsilon_3 + \frac{1}{\sqrt{2}}\epsilon_4)  + \epsilon_1\right) \Bigg)
\end{gather*}
-block-encoding of the coefficient matrix of system~\eqref{syst: normalized hermitian OSS} can be implemented efficiently according to Lemma 52 from \cite{gilyen2018quantum}.
We can choose
\begin{equation*}
    \begin{aligned}
    \epsilon_1 &= \frac{\epsilon_{QLSA} }{\kappa_M^3} \frac{1}{2\mathcal{K}}\\
    \epsilon_2 &= \frac{\epsilon_1}{2\sqrt{\|V\|_F^2 + \|A\|_F^2}}\\
    \epsilon_3 &= \epsilon_2\\
    \epsilon_4 &= \sqrt{2}\epsilon_2,
\end{aligned}
\end{equation*}
where $\mathcal{K}$ depends on the initial data
\begin{equation*}
    \mathcal{K} = \sqrt{2} \sqrt{\|V\|_F^2 + \|A\|_F^2}(\sqrt{2}\|Q\|_F+\sqrt{2}+1)^2.
\end{equation*}
Now, considering that the complexity for all the block-encoding algorithms we have used so far have poly-logarithmic dependence on the dimension and accuracy, and that, for $i=1,\ 2,\ 3,\ 4$ 
\begin{equation*}
    \begin{aligned}
    \mathcal{O}\left(\poly\log(\frac{1}{\epsilon_i})\right) &= \mathcal{O}\left(\poly\log(\kappa_M)\right),
\end{aligned}
\end{equation*}
the complexity for block-encoding will be dominated by the complexity for QLSA because QLSA has linear dependence on $\kappa_M$. So we can ignore the complexity of block-encoding.

\section[\appendixname~\thesection]{Spectral Analysis for Matrix $\Psi$}\label{sec: spectral analysis for the middle matrix}
In this section, we provide the spectral analysis for the matrix
\begin{equation}\label{eqn:middle matrix}
    \Psi = \begin{bmatrix}
        S^2 + 2\mu Q  & -\mu\theta E\\-\mu\theta E&X^2
    \end{bmatrix}.
\end{equation}
Just like in the previous section, for simplicity, we ignore the superscript $k$. We can do the following decomposition
\begin{equation*}
    \begin{aligned}
    \begin{bmatrix}
        S^2 + 2\mu Q  & -\mu\theta E\\-\mu\theta E&X^2
    \end{bmatrix} &= \begin{bmatrix}
        S^2   & -\mu\theta E\\-\mu\theta E&X^2
    \end{bmatrix} + \begin{bmatrix}
        2\mu Q  & 0\\0&0
    \end{bmatrix}.
\end{aligned}
\end{equation*}
Let us use the following notation
\begin{equation*}
    \begin{aligned}
    \Psi_1  &= \begin{bmatrix}
        S^2 & -\mu\theta E\\-\mu\theta E&X^2
    \end{bmatrix}\\
    \Psi_2  &= \begin{bmatrix}
        2\mu Q  & 0\\0&0
    \end{bmatrix}.
\end{aligned}
\end{equation*}
It can be proven that $\Psi_1$ is positive definite. The majority of the proof of this conclusion comes from the paper \cite{mohammadisiahroudi2021inexact}. For the reader's convenience, we provide the complete proof here.

Matrix $\Psi_1$ is a block diagonal matrix, with all the four blocks being diagonal matrices. So we can easily compute the eigenvalues using the characteristic polynomial
\begin{equation*}
\begin{aligned}
    \det(\Psi_1-q I) &= \det\Big(\big(X^2-q I\big)\big(S^2 -q I\big) - \theta^2\mu^2 E^2\Big)\\
                        &= \prod_{i=1}^n \Big(\big(x_i^2-q  \big)\big(s_i^2 -q \big) - \theta^2\mu^2 E_{ii}^2\Big).
\end{aligned}
\end{equation*}
Clearly, $\det(\Psi_1-q I)=0$ gives $n$ quadratic equations and each quadratic equation gives two eigenvalues. The two eigenvalues from the $i$th quadratic equation are
\begin{equation*}
    q_{i+} = \frac{1}{2}\Bigg((x_i^2 + s_i^2) + \sqrt{(x_i^2 + s_i^2)^2 - 4 x_i^2s_i^2 + 4\theta^2\mu^2E_{ii}^2}   \Bigg)
\end{equation*}
and
\begin{equation*}
     q_{i-} = \frac{1}{2}\Bigg((x_i^2 + s_i^2) - \sqrt{(x_i^2 + s_i^2)^2 - 4 x_i^2s_i^2 + 4\theta^2\mu^2E_{ii}^2}   \Bigg).
\end{equation*}
Recalling the definition of $E$ in \eqref{matrix E}, we can write
\begin{equation*}
\begin{aligned}
    q_{i-} &= \frac{1}{2}\Bigg((x_i^2 + s_i^2) - \sqrt{(x_i^2 + s_i^2)^2 - 4 x_i^2s_i^2 + 4(x_is_i - \mu)^2}   \Bigg)\\
    &= \frac{1}{2}\Bigg((x_i^2 + s_i^2) - \sqrt{(x_i^2 + s_i^2)^2 + 4(x_is_i - \mu + x_is_i)(x_is_i - \mu -x_is_i)}   \Bigg)\\
    &= \frac{1}{2}\Bigg((x_i^2 + s_i^2) - \sqrt{(x_i^2 + s_i^2)^2 - 4\mu(2x_is_i - \mu )}   \Bigg)\\
    &= \frac{1}{2}\Bigg((x_i^2 + s_i^2) - \sqrt{(x_i^2 + s_i^2)^2 - 4\mu(2\theta\mu E_{ii} + \mu )}   \Bigg).
\end{aligned}
\end{equation*}
One can verify that the square root always exists because
\begin{equation*}
\begin{aligned}
    (x_i^2 + s_i^2)^2 - 4\mu(2x_is_i - \mu )&\geq 4(x_is_i)^2 - 4 \mu (2x_is_i) + 4\mu^2\\
    &= 4(x_is_i-\mu)^2\\
    &\geq0.
\end{aligned}
\end{equation*}
With $\theta\in\left(0,\ \min \left\{\frac{1}{3\sqrt{n}},\ \frac{1}{4\|QVV^T\|_F + 1}\right\} \right)$, we have
\begin{equation*}
\begin{aligned}
    q_{i-} &\geq \frac{1}{2}\Bigg((x_i^2 + s_i^2) - \sqrt{(x_i^2 + s_i^2)^2 - 4\mu(2\theta\mu E_{ii} + \mu )}   \Bigg)\\
    &\geq \frac{1}{2}\Bigg((x_i^2 + s_i^2) - \sqrt{(x_i^2 + s_i^2)^2 - 4\mu(-2\mu\frac{1}{3\sqrt{n}}  + \mu )}   \Bigg)\\
    &= \frac{1}{2}\Bigg((x_i^2 + s_i^2) - \sqrt{(x_i^2 + s_i^2)^2 - \frac{4}{3}\mu^2}   \Bigg)\\
    &= \frac{1}{2}\frac{\frac{4}{3}\mu^2}{(x_i^2 + s_i^2) + \sqrt{(x_i^2 + s_i^2)^2 - \frac{4}{3}\mu^2} }\\
    &\geq \frac{1}{2}\frac{\frac{4}{3}\mu^2}{(x_i^2 + s_i^2) + \sqrt{(x_i^2 + s_i^2)^2} }\\
    &= \frac{\mu^2}{3(x_i^2+s_i^2)}\\
    &>0.
\end{aligned}
\end{equation*}
This means that matrix $\Psi_1$ is positive definite and its eigenvalues coincide with its singular values because $\Psi_1$ is also real and symmetric. Analogously, we have
\begin{equation*}
    \begin{aligned}
    q_{i+} &= \frac{1}{2}\Bigg((x_i^2 + s_i^2) + \sqrt{(x_i^2 + s_i^2)^2 - 4\mu(2\theta\mu E_{ii} + \mu )}   \Bigg)\\
    &\leq \frac{1}{2}\Bigg((x_i^2 + s_i^2) + (x_i^2 + s_i^2) + 2\mu\sqrt{(2\theta E_{ii} + 1 )}   \Bigg)\\
    &\leq \frac{1}{2}\Bigg((x_i^2 + s_i^2) + (x_i^2 + s_i^2) + 2\mu\sqrt{2}   \Bigg)\\
    &= (x_i^2 + s_i^2) + \sqrt{2}\mu. 
\end{aligned}
\end{equation*}
So the condition number of $\Psi$ satisfies
\begin{equation*}
    \begin{aligned}
    \kappa(\Psi) &\leq \frac{\sigma_{\max}(\Psi_1) + \sigma_{\max}(\Psi_2)}{\sigma_{\min}(\Psi_1) + \sigma_{\min}(\Psi_2)}\\
    &= \frac{\max_i q_{i+} + \sigma_{\max}(\Psi_2)}{\min_j q_{j-} + \sigma_{\min}(\Psi_2)}\\
    &\leq \frac{\max_i \{x_i^2+s_i^2\} +\sqrt{2}\mu + 2\mu\sigma_{\max}(Q)}{\min_{j} \frac{\mu^2}{3(x_i^2 + s_i^2)}}\\
    &= \frac{3\max_i \{x_i^2+s_i^2\} \left(\max_i \{x_i^2+s_i^2\} + \sqrt{2}\mu +2\mu\sigma_{\max}(Q) \right)}{\mu^2}\\
    &\leq \frac{3\omega^2 \left(\omega^2 +  \sqrt{2}\mu +2\mu\sigma_{\max}(Q)  \right)}{\mu^2},
\end{aligned}
\end{equation*}
where the last inequality comes from the definition of $\omega$. Since $\omega^2\geq x_is_i\geq(1-\theta)\mu$, we have
\begin{equation*}
    \kappa(\Psi) = \mathcal{O}\left(\frac{\omega^2(\omega^2 +  \mu\sigma_{\max}(Q))}{\mu^2}\right).
\end{equation*}
Using Lemma \ref{lemma:singular value QVQ}, we can also bound the condition number of matrix $M$ by
\begin{equation*}
    \begin{aligned}
    \kappa_M&=\sqrt{\kappa(M^TM)}\\
    &\leq \sqrt{\kappa(\Psi)}\kappa_{VAQ}\\
    &\leq \mathcal{O}\left(\frac{(\omega^2 +  \mu\sigma_{\max}(Q))}{\mu}\kappa_{VAQ}  \right).
\end{aligned}
\end{equation*}

\begin{adjustwidth}{-\extralength}{0cm}

\reftitle{References}


\bibliography{main}

\begin{thebibliography}{999}

\bibitem[Nocedal and Wright(1999)]{nocedal1999numerical}
Nocedal, J.; Wright, S.J.
\newblock {\em Numerical Optimization}; Springer,  1999.

\bibitem[Boser \em{et~al.}(1992)Boser, Guyon, and Vapnik]{boser1992training}
Boser, B.E.; Guyon, I.M.; Vapnik, V.N.
\newblock A training algorithm for optimal margin classifiers.
\newblock In Proceedings of the Fifth Annual Workshop on Computational Learning
  Theory; Haussler, D., Ed.,  1992, pp. 144--152.

\bibitem[Roos \em{et~al.}(1997)Roos, Terlaky, and Vial]{roos1997theory}
Roos, C.; Terlaky, T.; Vial, J.P.
\newblock {\em Theory and Algorithms for Linear Optimization: An Interior Point
  Approach}; John Wiley \& Sons,  1997.

\bibitem[P{\'o}lik and Terlaky(2010)]{polik2010interior}
P{\'o}lik, I.; Terlaky, T.
\newblock Interior point methods for nonlinear optimization. In {\em Nonlinear
  Optimization}; Gianni Di~Pillo, F.S., Ed.; Springer,  2010; pp. 215--276.

\bibitem[Gondzio(2013)]{gondzio2013convergence}
Gondzio, J.
\newblock Convergence analysis of an inexact feasible interior point method for
  convex quadratic programming.
\newblock {\em SIAM Journal on Optimization} {\bf 2013}, {\em 23},~1510--1527.

\bibitem[Lu \em{et~al.}(2006)Lu, Monteiro, and O'Neal]{lu2006iterative}
Lu, Z.; Monteiro, R.D.; O'Neal, J.W.
\newblock An iterative solver-based infeasible primal-dual path-following
  algorithm for convex quadratic programming.
\newblock {\em SIAM Journal on Optimization} {\bf 2006}, {\em 17},~287--310.

\bibitem[Bunch and Parlett(1971)]{bunch1971direct}
Bunch, J.R.; Parlett, B.N.
\newblock Direct methods for solving symmetric indefinite systems of linear
  equations.
\newblock {\em SIAM Journal on Numerical Analysis} {\bf 1971}, {\em
  8},~639--655.

\bibitem[Mohammadisiahroudi \em{et~al.}(2021)Mohammadisiahroudi, Fakhimi, and
  Terlaky]{mohammadisiahroudi2021inexact}
Mohammadisiahroudi, M.; Fakhimi, F.; Terlaky, T.
\newblock An Inexact Feasible Interior Point Method for Linear Optimization
  with High Adaptability to Quantum Computers.
\newblock {\em Technical Report} {\bf 2021}.

\bibitem[Harrow \em{et~al.}(2009)Harrow, Hassidim, and
  Lloyd]{harrow2009quantum}
Harrow, A.W.; Hassidim, A.; Lloyd, S.
\newblock Quantum algorithm for linear systems of equations.
\newblock {\em Physical Review Letters} {\bf 2009}, {\em 103},~150502.

\bibitem[Kerenidis and Prakash(2020)]{kerenidis2020quantum}
Kerenidis, I.; Prakash, A.
\newblock A quantum interior point method for LPs and SDPs.
\newblock {\em ACM Transactions on Quantum Computing} {\bf 2020}, {\em
  1},~1--32.

\bibitem[Mohammadisiahroudi \em{et~al.}(2022)Mohammadisiahroudi, Fakhimi, and
  Terlaky]{mohammadisiahroudi2022efficient}
Mohammadisiahroudi, M.; Fakhimi, R.; Terlaky, T.
\newblock Efficient Use of Quantum Linear System Algorithms in Interior Point
  Methods for Linear Optimization.
\newblock {\em arXiv preprint arXiv:2205.01220} {\bf 2022}.

\bibitem[Augustino \em{et~al.}(2021)Augustino, Nannicini, Terlaky, and
  Zuluaga]{augustino2021quantum}
Augustino, B.; Nannicini, G.; Terlaky, T.; Zuluaga, L.F.
\newblock Quantum interior point methods for semidefinite optimization.
\newblock {\em arXiv preprint arXiv:2112.06025} {\bf 2021}.

\bibitem[Schuld \em{et~al.}(2016)Schuld, Sinayskiy, and
  Petruccione]{SchuldRegression}
Schuld, M.; Sinayskiy, I.; Petruccione, F.
\newblock Prediction by linear regression on a quantum computer.
\newblock {\em Physical Review A} {\bf 2016}, {\em 94},~022342.

\bibitem[Kerenidis \em{et~al.}(2021)Kerenidis, Prakash, and
  Szil{\'a}gyi]{kerenidis2021quantum}
Kerenidis, I.; Prakash, A.; Szil{\'a}gyi, D.
\newblock Quantum algorithms for second-order cone programming and support
  vector machines.
\newblock {\em Quantum} {\bf 2021}, {\em 5},~427.

\bibitem[Kojima \em{et~al.}(1989)Kojima, Mizuno, and
  Yoshise]{kojima1989polynomial}
Kojima, M.; Mizuno, S.; Yoshise, A.
\newblock A polynomial-time algorithm for a class of linear complementarity
  problems.
\newblock {\em Mathematical Programming} {\bf 1989}, {\em 44},~1--26.

\bibitem[Monteiro and Adler(1989)]{monteiro1989interior}
Monteiro, R.D.; Adler, I.
\newblock Interior path following primal-dual algorithms. Part II: Convex
  quadratic programming.
\newblock {\em Mathematical Programming} {\bf 1989}, {\em 44},~43--66.

\bibitem[Goldfarb and Liu(1990)]{goldfarb1990n}
Goldfarb, D.; Liu, S.
\newblock An O ($n^3L$) primal interior point algorithm for convex quadratic
  programming.
\newblock {\em Mathematical programming} {\bf 1990}, {\em 49},~325--340.

\bibitem[Gily{\'e}n \em{et~al.}(2018)Gily{\'e}n, Su, Low, and
  Wiebe]{gilyen2018quantum}
Gily{\'e}n, A.; Su, Y.; Low, G.H.; Wiebe, N.
\newblock Quantum singular value transformation and beyond: exponential
  improvements for quantum matrix arithmetics [full version].
\newblock {\em arXiv preprint arXiv:1806.01838} {\bf 2018}, {\em 4}.

\bibitem[Chakraborty \em{et~al.}(2018)Chakraborty, Gily{\'e}n, and
  Jeffery]{chakraborty2018power}
Chakraborty, S.; Gily{\'e}n, A.; Jeffery, S.
\newblock The power of block-encoded matrix powers: improved regression
  techniques via faster Hamiltonian simulation.
\newblock {\em arXiv preprint arXiv:1804.01973} {\bf 2018}.

\bibitem[van Apeldoorn \em{et~al.}(2022)van Apeldoorn, Cornelissen, Gily{\'e}n,
  and Nannicini]{van2022quantum}
van Apeldoorn, J.; Cornelissen, A.; Gily{\'e}n, A.; Nannicini, G.
\newblock Quantum tomography using state-preparation unitaries.
\newblock {\em arXiv preprint arXiv:2207.08800} {\bf 2022}.

\bibitem[Horn and Johnson(2012)]{horn2012matrix}
Horn, R.A.; Johnson, C.R.
\newblock {\em Matrix analysis}; Cambridge university press,  2012.

\bibitem[Cortes and Vapnik(1995)]{cortes1995support}
Cortes, C.; Vapnik, V.
\newblock Support-vector networks.
\newblock {\em Machine Learning} {\bf 1995}, {\em 20},~273--297.

\bibitem[Rebentrost \em{et~al.}(2014)Rebentrost, Mohseni, and
  Lloyd]{rebentrost2014quantum}
Rebentrost, P.; Mohseni, M.; Lloyd, S.
\newblock Quantum support vector machine for big data classification.
\newblock {\em Physical Review Letters} {\bf 2014}, {\em 113},~130503.

\end{thebibliography}

\end{adjustwidth}
\end{document}